\documentclass[12pt,a4paper]{lms}

\usepackage{amsfonts}
\usepackage{amssymb}
\usepackage{diagrams}
\usepackage[pdftex,colorlinks=true,
                      pdfstartview=FitV,
                      linkcolor=blue,
                      citecolor=blue,
                      urlcolor=blue,
          ]{hyperref}

\newtheorem{thm}{Theorem}[section] 
\newtheorem{lemma}[thm]{Lemma}     
\newtheorem{cor}[thm]{Corollary}
\newtheorem{prp}[thm]{Proposition}

\newnumbered{assertion}{Assertion}    
\newnumbered{conjecture}{Conjecture}  
\newnumbered{lemma+df}{Lemma+Definition}
\newnumbered{df}{Definition}
\newnumbered{hypothesis}{Hypothesis}
\newnumbered{rem}{Remark}
\newnumbered{note}{Note}
\newnumbered{observation}{Observation}
\newnumbered{problem}{Problem}
\newnumbered{question}{Question}
\newnumbered{algorithm}{Algorithm}
\newnumbered{ex}{Example}
\newunnumbered{notation}{Notation} 

\newcommand{\mod}{{\rm mod}}
\newcommand{\Rad}{{\rm Rad}}
\newcommand{\Aut}{{\rm Aut}}

\newcommand{\Hom}{\mbox{{\rm Hom}}}

\newcommand{\res}{{\rm res}}

\newcommand{\B}{{\mathcal B}}

\newcommand{\F}{{\mathcal F}}

\renewcommand{\P}{{\mathcal P}}

\newcommand{\pf}{\smallskip {\bf Proof:\hspace{1em}}}
\newcommand{\epf}{$\circ$ \medskip}

\newcommand{\Soc}{{\rm Soc}}
\newcommand{\Dim}{{\rm Dim}}

\newcommand{\chr}{{\rm char}}

\title[Canonical Modules]
{Relative Invariants, Ideal Classes and Quasi-Canonical Modules of
Modular Rings of Invariants} 

\author{Peter Fleischmann and Chris Woodcock}

\classno{13B05, 20C20 (primary), 13B05,13B40,20C05 (secondary).
 Please refer to {\tt http://www.ams.org/msc/} for a list of codes}

\begin{document}
\maketitle

\begin{abstract}
We describe ``quasi canonical modules" for modular invariant
rings $R$ of finite group actions on factorial Gorenstein domains.
From this we derive a general ``quasi Gorenstein criterion" in terms
of certain $1$-cocycles. This generalizes a recent result
of A. Braun for linear group actions on
polynomial rings, which itself generalizes a classical
result of Watanabe for non-modular invariant rings.\\
We use an explicit classification of all reflexive rank one $R$-modules,
which is given in terms of the class group of $R$,
or in terms of $R$-semi-invariants. This result is
implicitly contained in a paper of Nakajima (\cite{Nakajima:rel_inv}).
\end{abstract}

\section{Introduction} \label{Intro}

Let $k$ be a field, $V$ a finite dimensional $k$-vector space of dimension $n$,
$G\subseteq {\rm GL}(V)$ a finite group and $A:={\rm Sym}(V^*)\cong k[x_1,\cdots,x_n]$,
the symmetric algebra over the dual space $V^*$ with its canonical $G$-action
and ring of invariants $R:=A^G:=\{a\in A\ |\ ga=a\ \forall g\in G\}.$\\
A classical result of K. Watanabe states that if $p=\chr(k)$ does not divide
$|G|$, then $A^G$ is Gorenstein if $G\subseteq {\rm SL}(V)$. If moreover
$G$ contains no pseudo-reflection, then the converse holds, i.e. if $A^G$ is Gorenstein,
then $G\subseteq {\rm SL}(V)$ (\cite{Watanabe:1}, \cite{Watanabe:2}).
In the recent paper \cite{Braun}, A. Braun proved an analogue of this result for the modular case, where
the characteristic of $k$ is allowed to divide the group order. Consider the following
\\[2mm]
{\bf Hypothesis} ($\mathcal{NR}$)\label{NR} :
The group $G\subseteq {\rm GL}(V)$ contains no pseudo-reflection (neither diagonalizable nor transvection).
\\[2mm]
Then Braun proved the following result:

\begin{thm}\cite{Braun}\label{Brauns_result}
Let $k$ be an arbitrary field and suppose the Hypothesis ($\mathcal{NR}$) holds.
Then the following are equivalent:
\begin{enumerate}
\item $G\subseteq {\rm SL}(V)$;
\item $A^G\cong {\rm Hom}_C(A^G,C)$ for every polynomial ring $C\subseteq A^G$ such that
$A^G$ is a finitely generated $C$-module and the homogeneous generators of $C$ have degrees divisible
by $|G|$.
\end{enumerate}
\end{thm}
From this he deduces that if $G$ satisfies Hypothesis ($\mathcal{NR}$), then the Cohen-Macaulay and Gorenstein loci
of $A^G$ coincide and if $A^G$ is Cohen-Macaulay it is also Gorenstein. He also obtains a modular
version of the converse: If $G$ satisfies Hypothesis ($\mathcal{NR}$) and $A^G$ is Gorenstein, then $G$ is
contained in ${\rm SL}(V)$.
\\[2mm]
In this paper we generalize Braun's results in two ways: firstly we avoid
Hypothesis ($\mathcal{NR}$) altogether. Secondly we neither assume $A$ to be a polynomial ring
nor that the parameter algebra $C$ is chosen in any particular way. Instead,
our main result applies, whenever $A$ is a (not necessarily graded) $k$-algebra, which
is also a factorial domain with unit group $U(A)=U(k)$.
It is remarkable that Braun's proof employs important techniques from the theory
of non-commutative Frobenius and symmetric algebras. The current paper grew out of
our attempt to understand these methods in detail, in particular with an eye on
possible future applications in non-commutative invariant theory. Nevertheless, it
turned out, that Braun's result as well as our generalization, can be obtained wholly within the
``world of commutative algebra", by combining Braun's ideas with information
hidden in the proofs of a classical paper by Nakajima \cite{Nakajima:rel_inv}.

To formulate our main result we need the following definitions and notation:

Let $A$ be a $k$-algebra which is also a factorial domain with
unit group $U(A)=U(k)$ and let $G\subseteq {\rm Aut}(A)$ be a finite group.
We do not assume that $G$ acts trivially on $k$, so $k':=k^G$ can
be a proper subfield of $k$.

\begin{df}\label{semiinv_df}
Let $\lambda\in Z^1(G,U(A))$ be a $1$-cocycle, i.e. $\lambda:\ G\to U(A)$ with
$$\lambda(gh)=\lambda(g)\cdot g(\lambda(h))\ \forall g,h\in G.$$
Then we define $A_\lambda:=\{a\in A\ |\ g(a):=\lambda(g)a\}$,
the \emph{$R$-module of relative $\lambda$-invariants}, or \emph{$\lambda$-semi invariants} in $A$.
\end{df}

\begin{df}\label{quasi_Gorenstein}
Let $\P$ be a commutative Gorenstein ring and $B$ a commutative
$\P$-algebra such that $_\P B$ is finite. Then we call
the $B$-module ${\rm Hom}_\P(B,\P)$ the \emph{quasi-canonical module}
of $B$ and we call $B$ \emph{quasi-Gorenstein} (w.r.t. $\P$), if ${\rm Hom}_\P(B,\P)\cong B$ as
$B$-modules (in other words, if $B$ is a symmetric $\P$-algebra).
\end{df}
\begin{rem} If $B$ is a graded connected $k$-algebra and
$\P$ a polynomial ring, generated by a homogeneous system of
parameters, then $B$ is a Cohen-Macaulay ring, if and only if
$_\P B$ is free. If $B$ is Cohen-Macaulay, then it is well known
that $\omega_B:={\rm Hom}_\P(B,\P)$ is a canonical module of $B$ and
$B$ is Gorenstein, if and only if $B\cong \omega_B$.
\end{rem}

Let $W:=W(G)\unlhd G$ be the normal subgroup generated by generalized reflections
(see Definition \ref{reflections on A}) and let $\F$ be any parameter $k'$-subalgebra $\F\subseteq R:=A^G\subseteq S:=A^W\subseteq A.$
Although not explicitly stated in \cite{Nakajima:rel_inv}, the following
facts are implicit in the proofs of that paper:
\\[2mm]
\begin{enumerate}\label{Nakajima_summary}
\item The class group ${\cal C}_R$ of $R$ is isomorphic to the subgroup
$\tilde H$ of $H^1(G,U(A))$, defined by
$\tilde H:=\{ \rho\in H^1(G,U(A))\ |\ {\rm res}_{_{I_{\rm Q}}}(\rho) = 1\
{\rm in}\ H^1(I_{\rm Q},U(A_{\rm Q})),\ \forall {\rm Q}\in {\rm Spec}_1(A)\}.$
(see Theorem \ref{faktoriell1}).
\item There are explicit bijections between the following sets:
\begin{itemize}
\item the divisor class group ${\cal C}_R$;
\item the set of iso types of finitely generated reflexive $R$-modules of rank one;
\item the set of iso types of $R$ modules of semi-invariants $A_\chi$ with
$\chi\in Z^1(G,U(A))$.
\end{itemize}
\item If $\chi\in Z^1(G/W,U(A))$, then $A_\chi\cong R \iff
[\chi]=1\in H^1(G/W,U(A))$.
\end{enumerate}

We can now state the main result of this paper:

\begin{thm}\label{main_intro}
The rings $S$ and $A$ are quasi-Gorenstein $\F$-algebras with
$${\rm Hom}_\F(S,\F)=S\cdot\theta_S\cong S,\ {\rm Hom}_\F(A,\F)=A\cdot\theta_A\cong A\ {\rm and}\
{\cal D}_{A,S}^{-1}\cong {\rm Hom}_S(A,S)= A\cdot \theta_{A,S}.$$
Here ${\cal D}_{A,R}={\cal D}_{A,S}$ is the Dedekind-different, which is a principal ideal in $A$ (see Definition \ref{Dedekind_diff}).
\\[1mm]
Let $\chi_S\in Z^1(G/W,U(k))$ and $\chi_A,\chi_{A,S}\in Z^1(G,U(k))$
be the ``eigen-characters" of $\theta_S$, $\theta_A$ and $\theta_{A,S}$, respectively.
Then $\chi_S=\chi_A\cdot\chi_{A,S}^{-1}$ and
${\rm Hom}_\F(R,\F)$ is isomorphic to the $R$-module of
semi-invariants $S_{\chi_S^{-1}}= A_{\chi_S^{-1}}$. In particular the following hold:
\begin{enumerate}
\item The quasi-canonical $R$-module ${\rm Hom}_\F(R,\F)$ is isomorphic to a divisorial ideal $I\unlhd R$, with
${\rm ch}({\rm cl}(I))=[\chi_S]=[\chi_A/\chi_{A,S}]$,
where ${\rm ch}:\ {\cal C}_R\to H^1(G/W,U(k))$
is the isomorphism of Corollary \ref{faktoriell2}.
\item The following are equivalent:
\begin{enumerate}
\item The ring $R$ is quasi-Gorenstein;
\item $[\chi_S]=1\in H^1(G/W,U(k))$.
\item $[\chi_{A,S}]=[\chi_A]\in H^1(G,U(k)).$
\end{enumerate}
\end{enumerate}
\end{thm}

\begin{rem}
\item In the special case, where $A$ is a polynomial ring with $k$-linear $G$ action,
the equivalence of (ii) (a) and (c) also appears in a paper by A Broer (\cite{Broer:TrGrps}).
\end{rem}

\begin{cor}\label{gor_cm_loc_intro}
If $[\chi_S]=1\in H^1(G/W,U(k))$, then the Cohen-Macaulay and Gorenstein loci of $R$ coincide.
\end{cor}

If $\chr(k)=p>0$, set $\tilde W:=\langle W,P^g\ |\ g\in G\rangle$ with $P$ a Sylow $p$-subgroup of $G$.
In other words, $\tilde W\unlhd G$ is the normal subgroup generated by all reflections on $A$ and all
elements of order a power of $p$. We obtain:

\begin{cor}\label{tilde_W}
If $G$ acts trivially on $k$, then $H^1(G/W,U(k))={\rm Hom}(G/W,U(k))={\rm Hom}(G/\tilde W,U(k))$
and Theorem \ref{main_intro} also holds with $W$ and $S$ replaced by $\tilde W$ and $\tilde S:=A^{\tilde W}$.
In particular $\tilde S$ is a factorial domain and quasi-Gorenstein and
$R$ is quasi-Gorenstein $\iff$ $\chi_{\tilde S}=1$.
\end{cor}

Assume for the moment that Hypothesis ($\mathcal{NR}$) holds, then $W=1$ and $A=S$ with
$[\chi_{A,S}]=1$. Hence in this case $R$ is quasi-Gorenstein, if and only
if $[\chi_A]=1$. If moreover $A={\rm Sym}(V^*)$ with $G\subseteq {\rm GL}_k(V)$, then
$[\chi_A]=\chi_A={\rm det}^{-1}$ (see Remark \ref{Sym_V_example}) and we recover
Braun's result (and Watanabe's for $\chr(k)\ \not |\ |G|$).
More generally:

\begin{cor}\label{refl_inv_poly_inv_gor}
Assume that $A={\rm Sym}(V^*)$ and $S:=A^W$ is Gorenstein (e.g a polynomial ring).
Assume moreover that $\chi_S=1$ (note that $\chi_S\in {\rm Hom}(G,U(k))$ here).
Then $R=A^G$ is Gorenstein, if it is Cohen-Macaulay.
\end{cor}

It is known by a result of Serre (\cite{bou}) that if
${\rm Sym}(V^*)^H$ is a polynomial ring for finite $H\subseteq {\rm GL}_k(V)$, then
$H=W(H)$. Unfortunately the converse is false, so the hypothesis of the above
Corollary is not automatic. If however it is satisfied, then the character
$\chi_S$ can be explicitly described in terms of the $G/W$ action on the homogeneous
generators of $A^W$ (see Section \ref{A_grd_conn_case}).

\section{The divisor class group and reflexive modules of rank one}
\label{div}

In this section we collect some definitions and results from \cite{Nakajima:rel_inv},
including some information which is implicitly contained via arguments
and proofs, but not explicitly stated there. In such a case we
include short proofs. Let $A$ be a Krull domain with quotient field $\mathbb{L}:= {\rm Quot}(A)$.
Let ${\rm Spec}_1(A):= \{ 0\ne {\rm P}\in {\rm Spec}(A)\ |\ {\rm ht(P)}=1\}$
then for every ${\rm P}\in {\rm Spec}_1(A)$, the localization $A_{\rm P}$ is
a discrete valuation ring and by definition
\begin{enumerate}
\item $A=\cap_{{\rm P}\in {\rm Spec}_1(A)} A_{\rm P}$;
\item for every $0\ne \ell\in \mathbb{L}$
the set $\{{\rm P}\in {\rm Spec}_1(A)\ |\ \nu_{\rm P}(\ell)\ne 0\}$ is finite.
\end{enumerate}
Let ${\cal D}_A$ denote the divisor group of $A$, i.e. the free abelian
group with basis ${\rm Spec}_1(A)$:
$${\cal D}_A:= \oplus_{{\rm P}\in {\rm Spec}_1(A)}\ \mathbb{Z}\ {\rm div}({\rm P}).$$
Let $0\ne J\triangleleft A$ be an ideal with $0\ne j\in J$. Then
$\nu_{\rm P}(J)\in \mathbb{Z}$ is defined by
$JA_{\rm P}={\rm P}^{\nu_{\rm P}(J)}A_{\rm P}$, hence
$\nu_{\rm P}(j):=\nu_{\rm P}(jA_{{\rm P}})\ge \nu_{\rm P}(J)\ge 0$,
and it follows that $\nu_{\rm Q}(J)=0$ for almost all ${\rm Q}\in {\rm Spec}_1(A)$.
If $I\subseteq \mathbb{L}$ is a fractional ideal, then
$\ell I\triangleleft A$ for some $\ell\in A$, hence again
$\nu_{\rm Q}(I)=0$ for almost all ${\rm Q}\in {\rm Spec}_1(A)$
and one defines
$${\rm div}(I):=\sum_{{\rm P}\in {\rm Spec}_1(A)} \nu_{\rm P}(I){\rm div}({\rm P}).$$
With ${\cal H}_A$ we denote the group of principal fractional ideals in $A$, then
the map ${\rm div}$ embeds ${\cal H}_A$ into ${\cal D}_A$ as a subgroup with
quotient group ${\cal C}_A:={\cal D}_A/{\cal H}_A$, the
\emph{divisor class group} of $A$.

\begin{df}\label{divisorial_closure}
Let $R\subseteq A$ be a subring.
For ideals $I\triangleleft R$ or $J\triangleleft A$ we denote with
$\overline I$ and $\overline J$ the corresponding divisorial closures, i.e.
$$\overline I=\cap_{Rr\triangleleft R\atop I\subseteq Rr}\ Rr,\ {\rm and}\
\overline J=\cap_{J\subseteq Aa}\ Aa.$$
\end{df}

\begin{lemma}\label{contract_divisorial}
Let $R\subseteq A$ be a subring with ${\rm Quot}(R)\cap A=R$, and let
$I\triangleleft R$ and $J\triangleleft A$ be ideals, then
\begin{enumerate}
\item[1.] $\overline{IA}\cap R=\overline I$;
\item[2.] $\overline J=d_JA$ with $d_J\in A$ $\iff$ $d_J:=\gcd(J)$ exists in $A$.
\item[3.] $\overline{IA}=d_IA$ with $d_I\in A$ $\iff$ $d_I:=\gcd(I):= \gcd \{r\in I\}\ ({\rm taken\ inside}\ A)$ exists.
\item[4.] For $a\in A$: $\overline{aJ}\subseteq \overline{J}\cdot a$ with equality,
if $J$ and $aJ$ have a $\gcd$ in $A$ (the latter is then $ad_J$).
\item[5.] For $a\in A$ and divisorial ideal $J\triangleleft A$,
$aJ\triangleleft A$ is divisorial.
\item[6.] If every subset of $A$ has a $\gcd$ (e.g. if $A$ is a factorial domain),
then for any ideals $J,K\triangleleft A$:
$\overline{\overline J\cdot K}=\overline{J\cdot K}$.
\end{enumerate}
\end{lemma}
\pf
By the assumption on ${\rm Quot}(R)$ we have $rA\cap R=rR$ for every $r\in R$. \\
1.:\ Let $\overline I=\cap_{I\subseteq rR\atop r\in R}\ rR$,
$\overline{IA}=\cap_{IA\subseteq aA\atop a\in A}aA$, and $x\in \overline{IA}\cap R$. Then
$I\subseteq rR$ implies $\overline{IA}\subseteq rA$, so $x\in rA\cap R=rR$.
It follows that $x\in \overline{I}$ and $\overline{I}\subseteq \overline{IA}\cap R\subseteq \overline{I}.$\\
2.+3.:\ Assume $d_J=\gcd(J)$; then clearly $J\subseteq Aa$ $\iff a\ |\ d_J$ $\iff d_JA \subseteq aA$.
Moreover, $J\subseteq d_JA$, so
$\overline J\subseteq d_JA\subseteq \cap_{J\subseteq Aa} Aa=\overline J.$
The opposite implications are obvious.\\
4. and 5.: $\overline{aJ}=$ $\cap_{Ja\subseteq Ac} Ac\subseteq$
$\cap_{Ja\subseteq Aba} Aba=$
$(\cap_{J\subseteq Ab} Ab)\cdot a=$ $\overline{J}\cdot a.$
If $J$ is divisorial we get
$\overline{aJ}\subseteq$ $\overline{J}\cdot a=aJ\subseteq \overline{aJ}.$
Let $g:=\gcd(J)$ and $d:=\gcd(aJ)$. Then $ag$ is a common divisor of $aJ$, hence
$ag$ divides $d$ and therefore $d/a\in A$ is a common divisor of $J$.
It follows that $d/a$ divides $g$, hence $d$ divides $ag$. So $ag | d | ag\sim d$.
Now we get $\overline{aJ}=\gcd(aJ)A=a\gcd(J)A=a\overline J.$\\
6.: For every $k\in K$ we have $\overline{J}k=\overline{Jk}\subseteq \overline{JK}$,
hence $\overline{J}K\subseteq \overline{JK}$ and
$\overline{\overline J\cdot K}\subseteq\overline{J\cdot K}$.
Clearly $JK\subseteq \overline{J}K$, hence
$\overline{JK}\subseteq \overline{\overline J\cdot K}.$
\epf

Let $B$ be an arbitrary commutative ring and $N\in B-mod$ a finitely generated
$B$-module. Then $N$ is torsion free of rank one $\iff$ there is an ideal $I\unlhd B$
containing a non zero-divisor, such that
$N\cong I\unlhd B$ are isomorphic as $B$-modules.

From now on let $A$ be a normal noetherian domain, then $A$ is a Krull-domain. Moreover
for every finitely generated module $M\in A-mod$ the following hold:
\begin{enumerate}
\item $M^*:={\rm Hom}_A(M,A)\cong \cap_{{\rm p}\in {\rm Spec}_1(A)}\ M_{\rm p}^* \subseteq \mathbb{L}\otimes_A M^*$.
\item If $M$ is torsion free, then the canonical map $c:\ M\to M^{**}$ induces an isomorphism
$$M^{**}\cong \cap_{{\rm p}\in {\rm Spec}_1(A)}\ M_{\rm p}.$$
\item The fractional ideal $I\in {\cal F}(A)$ is divisorial if and only if $I$ is a reflexive $A$ - module.
\item ${\rm ker}(c) = {\rm Tor}(M)$, the torsion submodule of $M$, and $M^*$ is reflexive.
\item For $M,N\in A-mod$ one has
$$({\rm Hom}_A(M,N))^{**}\cong {\rm Hom}_A(M^{**},N^{**}).$$
\end{enumerate}

\begin{prp}\label{classification_of_rank_one_reflexive}
Let $A$ be a normal noetherian domain, then there is a
bijection between the divisor class group
${\cal C}_A$ and the set of isomorphism classes
of finitely generated reflexive $A$-modules of rank one.
\end{prp}

\pf
If $M,N\in A-mod$ are f.g. reflexive $A$-modules of rank one, then
$M\cong I$ and $N\cong J$ with divisorial ideals $I,J\triangleleft A$, so we can
assume that $M=I,N=J$ are divisorial ideals.
Let $\theta:\ I\to J$ be an isomorphism, then for any $i,i'\in I$,
$\theta(ii')=i\theta(i')=i'\theta(i)$, so
$\ell:=\theta(i)/i\in \mathbb{L}$ with $\ell\cdot I\subseteq J$.
By symmetry we have
$\ell^{-1}=i/\theta(i)=$ $\theta^{-1}(\theta(i))/\theta(i)=$
$\theta^{-1}(j)/j$ for every $j\in J$, hence $j=\theta^{-1}(j))\ell$
and $J\subseteq \ell I$, so $J=\ell\cdot I$. It follows that the classes
${\rm cl}(J):=[{\rm div}(J)]$ and ${\rm cl}(I)\in {\cal C}_A$
coincide. \\
Now assume ${\rm cl}(J)={\rm cl}(I)\in {\cal C}_A$,
then ${\rm div}(I)={\rm div}(J)+{\rm div}(\ell A)$ for some
$\ell=a/b\in \mathbb{L}$, hence
$${\rm div}(Ib)={\rm div}(I)+{\rm div}(bA)={\rm div}(J)+{\rm div}(aA)=
{\rm div}(Ja),$$
and replacing $I$ by $Ib\cong I$ and $J$ by $Ja\cong J$, we can assume
that ${\rm div}(I)={\rm div}(J)$.
Hence $I_{\rm P}=J_{\rm P}$ for all ${\rm P}\in {\rm Spec}_1(A)$, so
$I\cong J$, since these are reflexive $A$-modules.
\epf

\section{Relative Invariants}\label{rel_inv}

Now let $G\subseteq \Aut(A)$ be a finite group of ring automorphisms with
corresponding ring of invariants $R:=A^G$ and quotient field
$\mathbb{K}=\mathbb{L}^G.$

The Galois group $G = {\rm Gal}(\mathbb{L}:\mathbb{L}^G)$ acts as permutation group on
${\rm Spec}_1(A)$ and on the
divisor group ${\cal D}_A$ and there is an inclusion homomorphism
$\rho:\ {\cal D}_{A^G}\to {\cal D}_A$
satisfying
$$d({\rm q}) \mapsto e_{\rm q}\cdot (\sum_{{\rm Q}\in {\rm Spec}_1(A):\ {\rm Q}\cap A^G={\rm q}}\ d({\rm Q})\ )
\in ({\cal D}_A)^G,$$
because the ramification index $e_{{\rm q},A}:=e_{\rm q}:= e({\rm Q}|{\rm q})$ is constant for all ${\rm Q}\in {\rm Spec}_1(A)$ over ${\rm q}$.
The group of invariants $({\cal D}_A)^G$ is a free abelian group with basis consisting of orbit sums
$$d({\rm Q})^+:= \sum_{\ g\in G / G_{\{\rm Q\}}}\ d(g{\rm Q}),\ {\rm Q}\in {\rm Spec}_1(A).$$
Here $G_{\{\rm Q\}}:={\rm Stab}_G({\rm Q})$ is the stabilizer (i.e. the decomposition group)
of ${\rm Q}$. Let $\mathfrak{C}$ denote a fixed set of representatives for the $G$-orbits on ${\rm Spec}_1(A)$,
i.e.
$$\mathfrak{C}\cong {\rm Spec}_1(A)/G\cong {\rm Spec}_1(A^G).$$
Then we have a short exact sequence of abelian groups:
\begin{equation}
\begin{diagram}
0 &\rTo&{\cal D}_{A^G}&\rTo^\rho&({\cal D}_A)^G&\rTo&\oplus_{{\rm Q}\in
\mathfrak{C}} {\mathbb{Z}}/e_{\rm q}{\mathbb{Z}}&\rTo&0 \\
\end{diagram}
\end{equation}

If $aA\in ({\cal H}_A)^G$, then $g(a)=c_g a$ with $c_g\in U(A)$ and
$gh(a)=c_{gh} a=$ $g(c_ha)=$ $g(c_h)c_g a$, hence
$c_{gh}=c_g\cdot g(c_h)$, so $\lambda:=c_{(\cdot)}\in Z^1(G,U(A))$ and $a\in A_\lambda$.

\begin{lemma}\label{Psi}
Let $\chi\in Z^1(G,U(A))$, then $0\ne A_\chi$ is a reflexive $R$-module of rank one and is
isomorphic to a divisorial ideal of $R$. The following hold:
\begin{enumerate}
\item for every $0\ne a\in A_{\chi^{-1}}$,
$\overline{aA_\chi A}\cap R=\overline{A_\chi a}=aA_\chi\triangleleft R$
is divisorial.
\item Let $\lambda\in B(G,U(A))$, i.e. $\lambda(g)=u^{-1}g(u)$ with $u\in U(A)$,
and $\mu:=\chi\cdot\lambda$. Then
$u\cdot A_\chi=A_\mu$ and $A_\mu A=A_\chi A,$ which only depends on the class $[\chi]\in H^1(G,U(A))$.
\item Assume $A$ to be a normal domain. Then
for every ${\rm Q}\in {\rm Spec}_1(A)$, $\nu_{\rm Q}(\overline{A_\chi A})<e({\rm Q}|{\rm q}).$
\end{enumerate}
\end{lemma}

\pf see \cite{Nakajima:rel_inv} Lemmas 2.1/2.2.
\epf

Let $Z_A^1(G,U(A)):=\{\lambda\in Z^1(G,U(A))\ |\ A_\lambda\not\subseteq {\rm Q}\ \forall {\rm Q}\in {\rm Spec}_1(A)\}$.
If $\lambda,\mu\in Z^1(G,U(A))$ and ${\rm Q}\in {\rm Spec}_1(A)$, then
$A_{\lambda\cdot\mu}\subseteq {\rm Q}$ implies
$A_{\lambda}\cdot A_{\mu}\subseteq A_{\lambda\cdot\mu}\subseteq {\rm Q}$,
hence $A_{\lambda}\subseteq {\rm Q}$, or
$A_{\mu}\subseteq {\rm Q}$. In other words, $Z_A^1(G,U(A))$ is a subgroup of $Z^1(G,U(A))$,
containing $B(G,U(A))$ (since $A_\lambda A=A$ for $\lambda\in B(G,U(A))$).
Therefore one can define
\begin{df}\label{H1_A}
$H_A^1(G,U(A)):=Z_A^1(G,U(A))/B(G,U(A)).$
\end{df}

\begin{lemma}\label{H1_sequence}
The sequence
\begin{diagram}\label{H1_sequence_diag}
0 &\rTo&H_A^1(G,U(A))&\rTo&H^1(G,U(A))&\rTo^{\Psi}&\oplus_{{\rm Q}\in
\mathfrak{C}} {\mathbb{Z}}/e_{\rm q}{\mathbb{Z}} \\
\end{diagram}
with $\Psi: [\chi]\mapsto (v_{\rm Q}(\overline{A_\chi A}))_{{\rm Q}\in
\mathfrak{C}}$
is an exact sequence of abelian groups.
\end{lemma}
\pf see \cite{Nakajima:rel_inv} Lemmas 2.3.
\epf

The map
$${\cal C}_{A^G} \to ({\cal D}_A)^G/({\cal H}_A)^G
\hookrightarrow ({\cal D}_A/{\cal H}_A)^G = ({\cal C}_A)^G$$
is essentially the natural map
$\phi:\ {\cal C}_{A^G} \to {\cal C}_A$ and we obtain
\begin{cor}\label{nak_2_4}
The kernel $\ker(\phi)$ is naturally isomorphic to $H_A^1(G,U(A))\cong\ker(\Psi)$.
Moreover, $\phi$ is
injective if and only if the $A_\chi$ are free $R$-modules for all
$\chi\in Z_A^1(G,U(A))$.
\end{cor}
\pf see \cite{Nakajima:rel_inv} Lemma 2.4.
\epf

\subsection{{\bf From now on we assume that $A$ is a factorial domain.}}

\begin{df}\label{reflections on A}
\begin{enumerate}
\item Let $I_{\rm Q}:=G_{k({\rm Q})}=\{g\in G\ |\ ga-a\in {\rm Q}\ \forall a\in A\}$, the inertia group of ${\rm Q}\in {\rm Spec}_1(A)$.
\item An element $g\in G$ is called a \emph{reflection on $A$}, if $g\in I_{\rm Q}$
for some ${\rm Q}\in {\rm Spec}_1(A)$. The group
$$W:=W_A:=W_A(G):=\langle I_{\rm Q}\ |\ {\rm Q}\in {\rm Spec}_1(A)\rangle$$
is a normal subgroup (since $G$ acts on ${\rm Spec}_1(A)$) and is called the \emph{subgroup of (generalized) reflections on $A$}.
\end{enumerate}
\end{df}

\begin{thm}\label{faktoriell1}
Let $A$ and $G$ be as above and assume that $A$ is a factorial domain. Then
${\cal C}_{A^G}\cong $
$$H_A^1(G,U(A)) \cong \tilde H:=\{ \rho\in H^1(G,U(A))\ |\ {\rm res}_{_{I_{\rm Q}}}(\rho) = 1\
{\rm in}\ H^1(I_{\rm Q},U(A_{\rm Q})),\ \forall {\rm Q}\in {\rm Spec}_1(A)\}.$$
In explicit form:\ Let $I$ be a divisorial ideal of $A^G$, then
$\overline{IA}=aA$ with semi-invariant $a\in A$. If $\theta_a\in Z^1(A,U(A))$
is the corresponding cocycle, i.e. $g(a)=\theta_a(g)a$ for every $g\in G$,
the class $[I]\in {\cal C}_{A^G}$ corresponds to the element $[\theta_a]\in \tilde H$.
\end{thm}
\pf See \cite{Nakajima:rel_inv} Lemma 2.4. The explicit form can be seen by following
up the isomorphism described there.
\epf

\begin{prp}\label{class_prp}
For $\chi\in Z^1(G,U(A))$ the following hold:
\begin{enumerate}
\item $\overline{A_\chi A}=d_\chi A$, $d_\chi:=\gcd(A_\chi)\in A_{\mu_\chi}$ with $\mu_\chi\in Z^1(G,U(A))$ and
a uniquely defined element $[\mu_\chi]\in H^1(G,U(A))$.
\item $A_\chi$ defines a unique class ${\rm cl}(A_\chi)\in {\cal C}_R$,
which satisfies
${\rm cl}(A_\chi)=[\chi^{-1}\mu_\chi]\in \tilde H$ (see \ref{faktoriell1}).
\item $A_\chi$ is a free $R$-module if and only if $[\chi]=[\mu_\chi]\in H^1(G,U(A))$.
\end{enumerate}
\end{prp}
\pf
(i):\ This follows from \ref{contract_divisorial}. \\
(ii):\ For every $a\in A_{\chi^{-1}}$ the ideal $aA_\chi\triangleleft R$ is divisorial
and we get from \ref{contract_divisorial}:\
$\overline{aA_\chi A}=a d_\chi A$ with $ad_\chi\in A_{\chi^{-1}\mu_\chi}.$
Hence ${\rm cl}(A_\chi)=[\chi^{-1}\mu_\chi]\in \tilde H$ by \ref{faktoriell1}.\\
(iii):\ This follows immediately from the above.
\epf

\begin{lemma}\label{cls_of_char2}
For $[\chi]\in H^1(G,U(A))$ the following are equivalent:
\begin{enumerate}
\item $[\chi]\in \tilde H=H^1_A(G,U(A))$;
\item $d_\chi\in U(A)$;
\item $[\chi^{-1}]={\rm cl}(A_{\chi})\in {\cal C}_R\cong \tilde H$;
\item $\overline{A_{\chi}A}=A.$
\end{enumerate}
\end{lemma}
\pf ``(i) $\iff$ (ii)":\
Let $[\chi]\in \tilde H$, then there is a divisorial ideal $J\unlhd R$ with
${\rm cl}(J)=[\chi^{-1}]$, i.e. $\overline{JA}=fA$ with $f\in A_{\chi^{-1}}$.
The divisorial ideal $I:= fA_\chi\triangleleft R$ satisfies
$$\overline{fA_\chi A}=\overline{IA}=f\cdot \overline{A_\chi A}=f d_\chi A.$$
Hence $J=\overline{JA}\cap R=$ $fA\cap R=fA_\chi=I$, so
$fd_\chi A=\overline{IA}=\overline{JA}=fA$ and $d_\chi\in U(A)$.
On the other hand, if $d_\chi\in U(A)$, then $[\mu_\chi]=1\in H^1(G,U(A))$ and
$[\chi^{-1}]=[\chi^{-1}][\mu_\chi]\in\tilde H$.\\
``(i) $\iff$ (iii)" and ``(ii) $\iff$ (iv)" follow from \ref{class_prp}.
\epf

\begin{cor}\label{princ_cor}
For $\chi\in Z^1(G,U(A))$ we have
$A_\chi=d_\chi\cdot A_{\chi\mu_\chi^{-1}}.$
Assume $A_\chi=a\cdot S$ with
$S\subseteq A$ and $a\in A$. Then $a\ |\ d_\chi$ and the following hold:
\begin{enumerate}
\item $a\sim d_\chi$ $\iff$ $S=A_\lambda$ with $[\lambda]=[{\chi\mu_\chi^{-1}}]\in \tilde H$ (i.e. $A_\lambda\cong A_1=R$ in $R-mod$.)
\item $1_A\in S$ $\iff$ $S=R$ $\iff$ $d_\chi\sim a\in A_\chi$.
\end{enumerate}
\end{cor}
\pf Since $d_\chi=\gcd(A_\chi)$, $A_\chi d_\chi^{-1}\subset A_{\chi\mu_\chi^{-1}}$, hence
$d_\chi\cdot A_{\chi\mu_\chi^{-1}}\subseteq A_\chi\subseteq d_\chi\cdot A_{\chi\mu_\chi^{-1}}$, so
$$A_\chi=d_\chi\cdot A_{\chi\mu_\chi^{-1}}.$$
(i):\ If $A_\chi=a S$ with $S\subseteq A\ni a$, then clearly $a\ |\ d_\chi$.
If $a=ud_\chi$ with $u\in U(A)$, then
$d_\chi\cdot A_{\chi\mu_\chi^{-1}}=A_\chi=$
$ud_\chi S$, hence $S=u^{-1}A_{\chi\mu_\chi^{-1}}=A_\lambda$ with
$[\lambda]=[{\chi\mu_\chi^{-1}}].$\\
Assume $S=A_\lambda\cong R$, then $d_\chi A=\overline{A_\chi A}=$
$\overline{aA_\lambda A}$=$a\overline{A_\lambda A}=aA$ by \ref{cls_of_char2}; hence
$a\sim d_\chi$.\\
(ii):\
If $1_A\in S$, then $a\in A_\chi$, therefore $d_\chi\ |\ a$ and
$S=1/a A_\chi\subseteq R$. Hence $A_\chi\subseteq aR\subseteq A_\chi$ and
$R=1/aA_\chi=S$.\\
If $S=R$, then $A_\chi=a R$, so $\gcd(A_\chi)\ni a \in A_\chi$.\\
If $d_\chi\sim a\in A_\chi$, $aR\subseteq A_\chi$, hence
$1_A\in R\subseteq 1/a A_\chi=S$.
\epf

\begin{cor}\label{dist_d}
Let $[\lambda]\in H^1(G,U(A))$ such that $A_\lambda=dR$. Then
for every $[\sigma]\in \tilde H$ we have
$$d=\gcd(A_{\lambda\sigma})\sim d_{\lambda\sigma},$$
i.e. $d$ and $d_{\lambda\sigma}$ are associated. In particular
$d=d_\lambda\cdot u$ with $u\in U(R)$ and $A_\lambda=Ad_\lambda.$

\end{cor}
\pf We have $dA=d\overline{A_\sigma A}=\overline{dRA_\sigma A}=
\overline{A_\lambda A_\sigma A}\subseteq \overline{A_{\lambda\sigma} A}=d_{\lambda\sigma}A=$
$$d_{\lambda\sigma}\overline{A_{\sigma^{-1}}A}=
\overline{d_{\lambda\sigma} A_{\sigma^{-1}}A}=
\overline{d_{\lambda\sigma}A A_{\sigma^{-1}}A}=\overline{\overline{A_{\lambda\sigma} A}A_{\sigma^{-1}}A}=\overline{A_{\lambda\sigma} A A_{\sigma^{-1}}A}\subseteq
\overline{A_\lambda A}=dA.$$
It follows that $d_\lambda=u\cdot d$ with $u\in U(A)\cap R=U(R)$.
\epf

\begin{cor}\label{dist_d}
Let $[\lambda]\in H^1(G,U(A))$ such that $A_\lambda=dR$. Then
for every $[\sigma]\in \tilde H$ we have
$$d=\gcd(A_{\lambda\sigma})\sim d_{\lambda\sigma},$$
i.e. $d$ and $d_{\lambda\sigma}$ are associated.
\end{cor}
\pf We have $dA=d\overline{A_\sigma A}=\overline{dRA_\sigma A}=
\overline{A_\lambda A_\sigma A}\subseteq
\overline{A_{\lambda\sigma} A}=d_{\lambda\sigma}A=d_{\lambda\sigma}\overline{A_{\sigma^{-1}}A}=$
$$\overline{d_{\lambda\sigma} A_{\sigma^{-1}}A}=
\overline{d_{\lambda\sigma}A A_{\sigma^{-1}}A}=
\overline{\overline{A_{\lambda\sigma} A}A_{\sigma^{-1}}A}=\overline{A_{\lambda\sigma} A A_{\sigma^{-1}}A}\subseteq
\overline{A_\lambda A}=dA.$$
\epf

\begin{thm}\label{canonical_cross_section}
Let $\P_{G,A}:=\{[\lambda]\in H^1(G,U(A))\ |\ A_\lambda=d_\lambda R\}$. Then
$$\P_{G,A}=\{[\lambda]\in H^1(G,U(A))\ |\ {\rm cl}(A_\lambda)=1\},\
\P_{G,A}\cap \tilde H=1\ {\rm and}\ H^1(G,U(A)) = \uplus_{[\lambda]\in \P_{G,A}} \tilde H\cdot [\lambda].$$
So $\P_{G,A}\subseteq H^1(G,U(A))$ is a transversal of the cosets of the subgroup $\tilde H\subseteq H^1(G,U(A))$. \\
For every $[\chi]\in H^1(G,U(A))$ let $[\mu_\chi]\in H^1(G,U(A))$ be the character of
$d_\chi:=\gcd(A_\chi)$, i.e. $d_\chi\in A_{\mu_\chi}$. Then the following hold
\begin{enumerate}
\item  ${\rm cl}(A_\chi)=[\chi^{-1}][\mu_\chi]$ with $\{[\mu_\chi]\}= \P_{G,A}\cap \tilde H\cdot [\chi]$.
\item The map $$\mu:\ H^1(G,U(A))\to \P_{G,A},\ [\chi]\mapsto [\mu_\chi]$$ satisfies
$\mu\circ\mu=\mu$ and it is a projection operator onto the distinguished
transversal $\P_{G,A}$.
\end{enumerate}
\end{thm}
\pf The equation $\P_{G,A}\cap \tilde H=1$  follows from
\ref{cls_of_char2} (iv).\\
Let $[\lambda],[\delta]\in \P_{G,A}$ with $[\sigma]:=[\lambda]^{-1}[\delta]\in \tilde H$, then
$[\delta]=[\lambda][\sigma]$, hence by Corollary \ref{dist_d},
$d_\delta\sim d_\lambda$ and $[\lambda]=[\delta]$. This shows that every $\tilde H$ coset contains at most one element
in $\P_{G,A}$.\\
Let $[\chi]\in H^1(G,U(A))$, then
$$\overline{A_\chi A}=d_\chi A\subseteq A_{\mu^{(1)}_\chi}A\subseteq \overline{A_{\mu^{(1)}_\chi}A}=
d_{\mu^{(1)}_\chi} A \subseteq  \overline{A_{\mu^{(2)}_\chi}A}=
d_{\mu^{(2)}_\chi} A \subseteq \overline{A_{\mu^{(3)}_\chi}A}=\cdots$$
with
$$[\chi]\equiv [\mu^{(1)}_\chi]\equiv [\mu^{(2)}_\chi]\equiv \cdots\ \mod\ \tilde H .$$
It is clear that this ascending chain of divisorial ideals must be stationary, hence
we will eventually have
$$d_{\mu^{(i)}_\chi}A=d_{\mu^{(i+1)}_\chi}A=d_{\mu^{(\infty)}_\chi},\ {\rm and\ therefore}\
[{\mu^{(i)}_\chi}]=[{\mu^{(i+1)}_\chi}]=[{\mu^{(\infty)}_\chi}]=:[\lambda]$$ with
$$\overline{A_\lambda A}=d_\lambda A \subseteq A_\lambda A\subseteq \overline{A_\lambda A}\ {\rm and}\
[\chi]\equiv [\mu_\chi]\equiv \cdots \equiv [\mu^{(\infty)}_\chi] =[\lambda]\ \mod\ \tilde H.$$
It follows that $d_\lambda=\gcd(A_\lambda)\in A_\lambda$, hence
$A_\lambda=d_\lambda R$, so $[\lambda]\in \P_{G,A}\cap \tilde H\cdot \chi.$\\
It now follows from Corollary \ref{dist_d} that
$$d_\lambda\sim d_{\mu^{(i)}_\chi}\sim d_{\mu^{(i-1)}_\chi}\sim d_{\mu^{(i-2)}_\chi}\sim \cdots\sim
d_{\mu^{(1)}_\chi}\sim d_\chi.$$
So $[\mu_\chi]:=[\mu^{(1)}_\chi]\in \P_{G,A}\cap \tilde H\cdot [\chi].$
By construction we have $d_{\mu_\chi}\sim d_\chi$, hence $\mu\circ \mu([\chi])=\mu([\chi])$,
which finishes the proof.
\epf

\begin{cor}\label{cross_section_of_ideal_classes}
For every $[\lambda]\in \P_{G,A}$ we have
${\cal C}_R = \{{\rm cl}(A_\chi)\ |\ \chi\in \tilde H\cdot [\lambda]\},$
i.e. if $\chi$ ranges through the full coset $\tilde H\cdot [\lambda]$, then the
$A_\chi$ form a transversal of all isomorphism types of rank one reflexive $R$-modules.\\
Alternatively the set $\{A_{\chi\mu_\chi^{-1}}\ |\ \chi\in Z^1(G,U(A))\}$ is also a full set of representatives of reflexive rank one $R$-modules.
\end{cor}
\pf
Every rank one reflexive $R$-module is isomorphic to a divisorial ideal of $R$,
the isomorphism type of which is uniquely determined by its ideal class.
From Corollary \ref{canonical_cross_section} we see that
$[\mu_{\sigma\lambda}]=$ ${\rm eigen character\ of}(d_{\sigma\lambda})=$
${\rm eigen character\ of}(d_\lambda)=[\lambda],$ hence we get
$${\rm cl}(A_{\sigma\lambda})=[\sigma]^{-1}[\lambda]^{-1}[\mu_{\sigma\lambda}]=
[\sigma]^{-1}[\lambda]^{-1}[\lambda]=[\sigma]^{-1}.$$
The last statement follows from \ref{princ_cor}, since
$A_\chi=d_\chi\cdot A_{\chi\mu_\chi^{-1}}\cong A_{\chi\mu_\chi^{-1}}$ in $R$-mod.
\epf

\subsection{{\bf $A$ noetherian, factorial domain, $U(A)=U(k)$}}\label{UFD_Uk}
{\bf From now on we assume that $A$ is a noetherian factorial domain with
$U(A)=U(k)$ with $k\subseteq A$, a field of characteristic $p\ge 0$. }
\\
Let ${\rm P}=a_{\rm P}A\in {\rm Spec}_1(A)$ and $\sigma\in I:=I_{\rm P}$.
Then for $u\in k$, $(\sigma-1)(u)\in k\cap {\rm P}=0$, so $\sigma(u)=u$ and
$W\subseteq {\rm Aut}_k(A)$. Clearly ${\rm P}$ is $I$-stable, so
$\sigma(a_{\rm P})=\delta_{\rm P}(g)a_{\rm P}$ and the map
$$\delta_{\rm P}:\ I_{\rm P}\to
U(k),\ \sigma\mapsto \delta_{\rm P}(g)=a_{\rm P}^{-1}\sigma(a_{\rm P})$$ is an element in $Z^1(I,U(k))={\rm Hom}(I,U(k))$.

\begin{lemma}\label{faktoriell2}
For ${\rm P}\in {\rm Spec}_1(A)$, $I:=I_{\rm P}$ and $e:=e({\rm P}|{\rm P}\cap R)$
we have ${\rm Hom}(I,U(k))={\rm Hom}(I,U(A_{\rm P}))=\langle\delta_{\rm P}\rangle\cong \mathbb{Z}/e \mathbb{Z}$.
There is a short exact sequence
\begin{diagram}\label{HUAS2}
0  &\rTo&{\cal C}_{A^G}&\rTo&H^1(G,U(k))&\rTo &\oplus_{{\rm Q}\in \mathfrak{C}}\
{\rm Hom}(I_{\rm Q},U(k))&\rTo&0
\end{diagram}
In particular ${\cal C}_{A^G} \cong H^1(G/W,U(k))$.
\end{lemma}
\pf see \cite{Nakajima:rel_inv} Lemmas 2.6. In addition to this, we only
need to show that $\tilde H=H^1(G/W,U(k)).$
Let $[\chi]\in\tilde H$ with $\chi\in Z^1(G,U(k))$, then for $g,h\in W$,
$\chi(gh)=\chi(g)g(\chi(h))=\chi(g)\chi(h)$, since $W$ acts trivially on $k$.
Moreover $g$ and $h$ are products of elements on which $\chi$ is $1$, hence
$\chi_{|W}=1$. We view $Z^1(G/W,U(k))$ as a subset of $Z^1(G,U(k))$ in a
natural way. Then, again since $W$ acts trivially on $k$ we have
$B^1(G,U(k))\subseteq Z^1(G/W,U(k))$, hence $B^1(G,U(k))=B^1(G/W,U(k))$, so
${\cal C}_{A^G}=\tilde H\cong Z^1(G/W,U(k))/B^1(G/W,U(k))=H^1(G/W,U(k)).$
\epf

\subsection{{\bf $A$ noetherian, factorial domain, $U(A)=U(k)$ with trivial $G$-action}}

Then $H^1(G,U(A))=G^*:={\rm Hom}(G,U(k))$, the group of linear
$k$-characters of $G$. If $N\unlhd G$ is a normal subgroup, then the restriction map yields a short
exact sequence
$$1\to (G/N)^*\to G^* \to N^*\to 1.$$
\begin{cor}\label{grad_con_ufd}
There is an isomorphism
${\rm ch}:\ {\cal C}_{A^G} \cong {\overline G}^* = \ker(\res_{|W}),$
where $\overline G:=G/W$ and $\res_{|W}:\ G^*\to W^*$ is
the restriction map on characters.
\end{cor}

\section{Quasi-Gorenstein Rings of Invariants}\label{q_gor_inv}

Now let $k$ be a field and $A$ a finitely generated normal $k$-algebra with
$U(A)=U(k)$, such that the quotient field $\mathbb{L}:={\rm Quot}(A)$ is separable over $k$.
Let $G\subseteq {\rm Aut}(A)$ be a finite group with ring of invariants
$R:=A^G$. Then $k$ is a separable algebraic extension of
$k':=k^G\subseteq \mathbb{K}:={\rm Quot}(R)$ and $\mathbb{L}$ as well as $\mathbb{K}$
are separable over $k'$. By Noether-normalization
there is a $k'$-polynomial ring $\F\subseteq R:=A^G$ such that ${_\F}R$ and $_{\F}A$ are
finitely generated modules, i.e. $\F=k'[f_1,\cdots,f_d]$, with $(f_1,\cdots,f_d)$ a system of parameters of $R$ as $k'$-algebra. It follows from \cite{eis} Cor. 16.18 pg.403, that
$\F$ can be chosen such that $\mathbb{L}$ and $\mathbb{K}$ are separable over
${\rm Quot}(\F)$. For technical reasons, which become clear later in section \ref{Dedekind_diff},
we choose and fix $\F$ in such a way.

\begin{prp}\label{cm_loc_is_gor_loc}
Let $\P\subseteq A$ as above be quasi-Gorenstein. Then for every
${\rm Q}\in {\rm Spec}(A)$, the localisation
$A_{\rm Q}$ is Cohen-Macaulay if and only if $A_{\rm Q}$ is Gorenstein.
In other words, the Cohen-Macaulay and Gorenstein loci of $A$ coincide.
\end{prp}
\pf
Let ${\rm Q}\in {\rm Spec}(A)$ be such that $A_{\rm Q}$ is Cohen-Macaulay.
Set ${\rm q}={\rm Q}\cap \P\in {\rm Spec}(\P)$ and let
${\rm Q}:={\rm Q}_1,\cdots,{\rm Q}_k$ be the primes of $A$ lying over ${\rm q}$.
Since ${\rm Hom}_\P(A,\P)\cong A$ and
$\widehat{A_{\rm q}}\cong \times_{i=1}^k \hat A_{{\rm Q}_i}$,
we get
$\widehat{A_{\rm q}}\cong \widehat{({\rm Hom}_\P(A,\P))_{\rm q})}\cong$
${\rm Hom}_\P(A,\P)\otimes_\P \hat \P_{\rm q}\cong$
${\rm Hom}_{\hat \P_{\rm q}}(\widehat{A_{\rm q}},\hat \P_{\rm q})\cong$
$\times_{i=1}^k {\rm Hom}_{\hat \P_{\rm q}}(\hat A_{{\rm Q}_i},\hat \P_{\rm q})\cong$
$\times_{i=1}^k \hat A_{{\rm Q}_i}$.
Let $1_{\widehat{A_{\rm q}}}=\sum_{i=1}^k e_i$
with $e_ie_j=\delta_{ij}$, then
$\hat A_{{\rm Q}_i}=e_i\widehat{A_{\rm q}}\cong$
$e_i{\rm Hom}_{\hat \P_{\rm q}}(\widehat{A_{\rm q}},\hat \P_{\rm q})\cong$
${\rm Hom}_{\hat \P_{\rm q}}(e_i\widehat{A_{\rm q}},\hat \P_{\rm q})\cong$
${\rm Hom}_{\hat \P_{\rm q}}(\hat A_{{\rm Q}_i},\hat \P_{\rm q})$.
Since $A_{{\rm Q}_1}$ is Cohen-Macaulay, so is
$\hat A_{{\rm Q}_1}$ and it is finite over $\hat \P_{\rm q}$.
It follows that
${\rm Hom}_{\hat \P_{\rm q}}(\hat A_{{\rm Q}_i},\hat \P_{\rm q})$
is the unique canonical module $\omega_{\hat A_{{\rm Q}_i}}$ (up to isomorphism) of
$\hat A_{{\rm Q}_i}$. Therefore
$\omega_{\hat A_{{\rm Q}_i}}\cong \hat A_{{\rm Q}_i}.$
It is generally true, that for a finitely generated $A_{\rm Q}$-module
$M$, the completion $M\otimes_{A_{\rm Q}} \hat A_{\rm Q}$ is canonical
for $\hat A_{\rm Q}$, if and only if $M$ is canonical for $A_{\rm Q}$,
so we conclude that $\omega_{A_{\rm Q}}\cong A_{\rm Q}$ and
$A_{\rm Q}$ is Gorenstein.
\epf

\begin{df}\label{Dedekind_diff}
For a normal subring $S\subseteq A$ such that $S\hookrightarrow A$ is finite
and ${\rm Quot}(A)$ is separable over ${\rm Quot}(S)$ let
${\mathcal D}_{A,S}\unlhd A$ denote the corresponding Dedekind different.
\end{df}

It is well known that ${\mathcal D}_{A,S}$ and its inverse
${\mathcal D}_{A,S}^{-1}$ are divisorial (fractional) ideals with
$${\mathcal D}_{A,S}^{-1}\cong {\rm Hom}_S(A,S).$$

{\bf Now we assume in addition that $A$ is a factorial domain (see subsection \ref{UFD_Uk})}.\\
Let $S:=A^W$, then by \ref{faktoriell2} $S$ is also factorial.
The following lemma is well known (at least in the context
of Dedekind domains appearing in number theory):

\begin{lemma}\label{gen_reflection_subgroup}
For any $W\subseteq H\subseteq G$ the following holds:
\begin{enumerate}
\item ${\mathcal D}_{A,A^G}={\mathcal D}_{A,A^H}.$
\item ${\mathcal D}_{A^H,A^G}=(1)=A^H.$
\end{enumerate}
In particular the extension $A^G\hookrightarrow A^W$ is unramified
in height one.
\end{lemma}

Using the fact that $_RS$ is unramified in height one we can now prove
the main result:

{\bf Proof} (of Theorem \ref{main_intro}):  We have ${\cal D}_{S,R}=S$,
$R^*={\cal D}_{R,\F}^{-1}$ and ${\cal D}_{S,\F}=Sd,$ a principal ideal,
since $S$ is a factorial domain. It follows that
$S^*=S\theta_S$, where $\theta_S\in S^*$
can be identified with an element in ${\rm Quot}(S)$.
Since the fractional ideal ${\cal D}_{S,\F}^{-1}$ is $G/W$-stable
$\theta_S$ is a relative invariant with character
$\chi_S\in\tilde H$.
By the Dedekind-tower theorem,
${\cal D}_{S,\F}=\overline{{\cal D}_{S,R}{\cal D}_{R,\F}}\unlhd S$,
which implies (first locally at height one primes, then globally):
$$S^*\cong S\theta_S={\cal D}_{S,\F}^{-1}=\overline{{\cal D}_{S,R}^{-1}S{\cal D}_{R,\F}^{-1}}=
\overline{S{\cal D}_{R,\F}^{-1}}\subseteq {\rm Quot}(S).$$
There is a suitable element $r\in R$ with $rS\theta_S\subseteq S$ and therefore
$rS\theta_S=r\overline{S{\cal D}_{R,\F}^{-1}}\subseteq S.$
Hence we get
$r\overline{S{\cal D}_{R,\F}^{-1}}\cap R=r{\cal D}_{R,\F}^{-1}=rS\theta_S\cap R$,
so $R^*\cong {\cal D}_{R,\F}^{-1}=S\theta_S\cap {\rm Quot}(R)=S_{{\chi_S}^{-1}}=A_{{\chi_S}^{-1}},$
where the isomorphism is one of $R$-modules.
Since $\chi_S\in \tilde H$ we have $[\mu_{\chi_S}]=1$, so
${\rm ch}({\rm cl}(A_{{\chi_S}^{-1}}))=[\chi_S].$
The equation $\chi_S=\chi_A\cdot\chi_{A,S}^{-1}$ follows immediately from
$$\overline{{\cal D}_{A,R}{\cal D}_{R,\F}}=
\overline{{\cal D}_{A,S}{\cal D}_{S,\F}}$$
and ${\cal D}_{A,S}={\cal D}_{A,R}.$
The remaining statements follow immediately.
\epf

{\bf Proof} of Corollary \ref{tilde_W}:
Since $\tilde W/W$ is generated by $p$-elements, it follows that
${\cal C}_{\tilde S}={\rm Hom}(\tilde W/W,U(k))=1$, hence
${\rm Hom}(G/W,U(k))={\rm Hom}(G/\tilde W,U(k))$ and
$\tilde S$ is a factorial domain, hence quasi Gorenstein. Using Lemma \ref{gen_reflection_subgroup}
the remaining arguments are exactly as above with $W$ replaced by $\tilde W$ and
$S$ by $\tilde S$.
\epf

{\bf Proof} of Corollary \ref{gor_cm_loc_intro}:
this follows immediately from Theorem \ref{main_intro} and
Proposition \ref{cm_loc_is_gor_loc}.
\epf

\section{The graded connected case}\label{A_grd_conn_case}

The application of Theorem \ref{main_intro} depends on the determination
of $[\chi_S]$ or, equivalently $[\chi_A]$ and $[\chi_{A,S}]$.
If $G$ acts trivially on $k$, then these are linear characters in
${\rm Hom}(G/W,U(k))$ or ${\rm Hom}(G,U(k))$, respectively. In this section
we investigate these characters in the case where $A$ is a graded connected Cohen-Macaulay
ring. \\
So throughout this section $A=\sum_{i\ge 0} A_i$ is an $\mathbb{N}_0$
graded connected noetherian normal $k$-algebra, i.e. $A_0=k$ with $U(A)=U(k)$
and $G\subseteq \Aut_k(A)$ a finite group of graded $k$-algebra automorphisms.
We will also assume that $A$ is a Cohen-Macaulay domain, i.e.
$A$ is a free module over some (and then every) parameter algebra $\F\subseteq A$.
We keep the previous notation, so
$R=A^G\hookrightarrow A$ is a finite extension of noetherian normal domains.
Let $y_1,y_2,\cdots,y_d\in R$ be a homogeneous system of parameters (hsop)
with $d_i:=\deg(y_i)$, $d=\Dim(R)=\Dim(A)$, and set $\F:=k[y_1,\cdots,y_d]$.

\begin{df}\label{hilb_tr}
Let $V:=\oplus_{n\ge 0} V_n$ be an $\mathbb{N}_0$ graded $k$-vectorspace
and $G$ a finite group acting on $V$ by graded $k$-linear automorphisms.
We define the \emph{(Brauer-) character series}
$$H_{V,g}^{(Br)}(t):=\sum_{n=0}^\infty \chi_{V_n}(g)t^n,$$
where $\chi_{V_n}$ is the (Brauer-) character afforded by the action of $G$ on $V_n$.
Note that $H_{V,g}(t)\in k[[t]]$, whereas $H_{V,g}^{Br}(t)\in \mathbb{Q}(\epsilon)[[t]]$,
where $\epsilon$ is a primitive ${\rm order}(g)$-th root of unity in $\mathbb{C}$.
\end{df}

Note that
$$H_A(t):=H^{Br}_{A,{\rm id}}(t)=\sum_{i\ge 0} \dim_k(A_i)t^i\in \mathbb{Q}(t)$$
is the ordinary Hilbert-series of $A$. Let
$$U:=\overline A:= A/\F^+A,$$ where
$\F^+:=(y_1,\cdots,y_d)\unlhd \F$ is the unique maximal homogeneous
ideal of $\F$. Then $\F\otimes_k U$ is the projective cover of $_\F A$ in $\F-mod$,
hence, as $_\F A$ is free, we have $\F\otimes_k U\cong A$ as $\F$-modules.
Moreover
$$U=\oplus_{i=0}^\beta U_i = \oplus_{i=1}^\ell k\xi_i,$$
where we choose a homogeneous $k$-basis $\{\xi_i\ |\ i=1,\cdots,\ell\}$
with $\deg(\xi_i)=:\beta_i\le \beta_{i+1}$, $\beta:=\beta_\ell$ and
$\ell:=\dim_k(U)$. We also will choose an $\F$-basis $\B:=\{s_i\ |\ i=1,\cdots,\ell\}$
of $A$, such that $s_i+\F^+A=\overline{s_i}=\xi_i$ for $i=1,\cdots,\ell$.\\
Note that $G$ acts on $A$ and $U$ and if
$g(\xi_i)=\sum_{j=1}^\ell g_{ji}\xi_j$ with $(g_{ji})\in k^{\ell\times \ell}$,
then
$$g(s_i)=\sum_{j=1}^\ell g_{ji}s_j + {\mathcal X}$$
with ${\mathcal X}\in \F^+A$. For each $j$ let $\tilde A_j:=\langle \B\rangle_k\cap A_j$,
then
$A_i=\oplus_{m+n=i} \F_m\otimes_k \tilde A_n$ and it is easily seen that
$$\chi_{A_i}(g)=\sum_{m+n=i} \dim_k(\F_m)\cdot \chi_{U_n}(g)=
{\rm coeff}_i(H^{Br}_\F(t)\cdot H^{Br}_{U,g}(t)).$$
Hence
$$H^{Br}_{A,g}(t)=H_{\F}(t)\cdot H^{Br}_{U,g}(t).$$
Since $H_\F(t)=\frac{1}{\prod_{i=1}^\ell(1-t^{d_i})}$ and $H^{Br}_{U,g}(t)\in \mathbb{Q}(\epsilon)[t]$, we get
\begin{lemma}\label{rat_br_hilb_ser}
The Brauer-character series of $A$ are rational, i.e.
$H^{Br}_{A,g}(t)\in \mathbb{Q}(\epsilon)(t)$.
\end{lemma}

{\bf Now we assume in addition that $A$ is Gorenstein.} It is then well known that
$$H_A(t)=(-1)^dt^{a(A)}H_A(1/t),$$
where $a(A)=\deg(H_A(t))$ is the degree of $H_A(t)$.
This symmetry is induced by the duality of the corresponding
artinian Gorenstein algebra
$$U=\overline A:= A/\F^+A,$$ where
$\F^+:=(y_1,\cdots,y_d)\unlhd \F$ is the unique maximal homogeneous
ideal of $\F$. For later use we recall the details:
\\[2mm]
There is a graded embedding
$$U/U^+[-\beta]\hookrightarrow k[-\beta]\subseteq\ _UU_\beta,\
k=((U/U^+)[-\beta])_\beta\ni \lambda\mapsto \lambda \xi_\ell.$$
It follows from \cite{BH} that $_UU$ is injective with
$\Soc(U)\cong k$ (up to shift), hence
$$k[-\beta]\cong U/U^+[-\beta]\cong {\rm Soc}(_UU).$$
It is well known that $^*E(k) \cong U^*:= {\rm Hom}_k(U,k)$,
where $^*E(k)$ denotes the graded $^*$injective hull of $k=U_0$
(see \cite{BH} for the definition of $^*$injectivity).
Note that $U=\oplus_{i=0}^\beta U_i$; choosing a homogeneous dual
$k$-basis $\{\xi_i^*\ |\ i=1,\cdots,\ell\}$ (such that $\xi_i^*(\xi_j)=\delta_{i,j}$ and
$\deg(\xi_i^*)=-\deg(\xi_i)$), we see that
$U^*=\oplus_{i=0}^\beta (U^*)_{-i}$ with $k\cong {\rm Soc}(U^*)\cong U^*_0$
and ${\rm dim}_k (U^*)_{-i}={\rm dim}_k U_i$.
Since $_UU$ is injective and indecomposable we conclude
$$_UU\cong\ ^*E({\rm Soc}(_UU))\cong
\ ^*E(k[-\beta])\cong\ ^*E(k)[-\beta]\cong U^*[-\beta].$$
It follows that
${\rm dim}_k(U_i)=$ ${\rm dim}_k(U^*[-\beta]_i)=$
${\rm dim}_k((U^*)_{i-\beta})=$ ${\rm dim}_k(U_{\beta-i})$,
hence $H_U(t)=t^\beta H_U(1/t)=H_U^*(t)$.
Since $\Rad(U^*)=\Soc(U)^\perp=\langle \xi_1^*,\cdots,\xi_{\ell-1}^*\rangle$
we have $_UU^*= U\cdot \xi_\ell^*$ as well as a non-degenerate
associative bilinear form
$$\kappa(\ ,\ ):\ U\times U\to k,\ \kappa(\xi,\xi')=\xi_\ell^*(\xi\cdot \xi').$$
It follows from $\Soc(U)=k\cdot\xi_\ell$, that for $g\in G$, $g(\xi_\ell)=\lambda(g)\xi_\ell$,
with some linear character $\lambda\in \Hom(G,U(k))$.
Since the $G$-action preserves degrees,
we have $g(\xi_j)\in \sum_{n<\beta} U_n$, hence
$g^{-1}\xi_\ell^*(\xi_j)=\xi_\ell^*(g(\xi_j))=0$
for every $j<\ell$ and
$$g^{-1}\xi_\ell^*(\xi_\ell)=\xi_\ell^*(g(\xi_\ell))=\lambda(g)\cdot 1;$$
hence $g\xi_\ell^*=\lambda(g)^{-1}\xi_\ell^*$ for every $g\in G$.
It follows that
$$\kappa(g(\xi_i),g(\xi_j))=\lambda(g)\cdot \kappa(\xi_i,\xi_j).$$

\begin{prp}\label{hilb_tr_quot} Let $A$ be a graded connected Gorenstein algebra, then
the Brauer-character series of $A$ and $U$ satisfy the following identities:
\begin{enumerate}
\item $H^{Br}_{U,g}(t)=\hat\lambda(g)\cdot t^\beta H^{Br}_{U,g^{-1}}(1/t)$;
\item $\frac{H^{Br}_{A,g}(t)}{H^{Br}_{A,g^{-1}}(1/t)}=(-1)^d t^{a(A)}\hat\lambda(g)$
with $a(A)=\beta-\sum_i d_i=\deg(H^{Br}_{A,1}(t)).$
\end{enumerate}
In particular
$$\hat\lambda(g)=(-1)^d\cdot \lim_{t\to 1}\frac{H^{Br}_{A,g}(t)}{H^{Br}_{A,g^{-1}}(1/t)}.$$
\end{prp}
\begin{rem}\label{lambda_A}
It follows from (i) that the character $\lambda$ only depends on $A$ and not on the
choice of $\F$. Therefore we denote it by $\lambda_A$ and we will denote the corresponding
Brauer character by $\hat\lambda_A$.
\end{rem}
\pf
1.:\ Let ${\frak A}:=\{a_1,\cdots,a_m\}$ and ${\frak B}:=\{b_1,\cdots,b_m\}$
be $k$-bases of $U_i$ and $U_{\beta-i}$, respectively, then
$\kappa(g(a_i),g(b_j))=\lambda(g)\cdot \kappa(a_i,b_j).$ On the other hand,
this is equal to
$M_{\frak A}(g)^{tr}\circ Q\circ M_{\frak B}(g)$, where
$Q=\kappa(a_i,b_j)\in k^{m\times m}$. For every $0\le \nu\le\beta$ with
$\nu\ne \beta-i$ we have
$$U_\nu\subseteq U_i^\perp:=\{a\in U\ |\ \kappa(a,U_i)=0\}.$$
Hence the map
$$U_i\times U_{\beta-i}\to k,\ (a,b)\mapsto \kappa(a,b)$$
is a perfect pairing, in particular $Q$ is a non-singular matrix.
Therefore
$M_{\frak A}(g)^{tr}=\lambda(g)\cdot Q\circ M_{\frak B}(g)^{-1}\circ Q^{-1}$
and $${\rm trace}(g_{|U_i})=\lambda(g) {\rm trace}(g^{-1}_{|U_{\beta-i}}),$$
from which 1. follows immediately.\\
2.:\ Using 1., the LHS is equal to
$$\frac{H^{Br}_{U,g}(t)}{H^{Br}_{U,g^{-1}}(1/t)}\cdot
\frac{\prod_i (1-t^{d_i})^{-1}}{\prod_i (1-t^{-d_i})^{-1}}=
\hat\lambda(g)\cdot t^{\beta-\sum_i d_i}(-1)^d.$$
\epf

\begin{rem}\label{Sym_V_example}
Let $g\in {\rm GL}(V)$ semisimple, $A:={\rm Sym}(V^*)\cong k[x_1,\cdots,x_n]$ with
$x_1,\cdots,x_n$ a basis of $V^*$. We can assume that $g(x_i)=\lambda_i x_i$ with
eigenvalues $\lambda_i\in U(k)$, so with slight abuse of notation we obtain
$$H^{Br}_{A,g}(t)=\widehat{{\rm trace}(g_{|A})}= \prod_{i=1}^n (1+\widehat{\lambda_i}t+\widehat{\lambda_i}^2t^2+\cdots)=
\prod_{i=1}^n \frac{1}{1-\widehat{\lambda_i} t}=\frac{1}{\widehat{\det(1-tg)}}.$$
It follows that
$H^{Br}_{A,g^{-1}}(1/t)=$ $\frac{1}{\widehat{\det(1-g^{-1}1/t)}}=$
$$\frac{t^n}{\widehat{\det(t- g^{-1})}}=\frac{t^n\widehat{\det(g)}}{\widehat{\det(gt-1)}}=
(-1)^nt^n\widehat{\det}(g)\cdot \frac{1}{\widehat{\det(1-gt)}}=
(-1)^nt^n\widehat{\det}(g)\cdot H_{A,g}(t).$$
Hence $\widehat{\lambda_A}(g)=\widehat{\det}(g)^{-1}.$
\end{rem}

\begin{prp}\label{coprime}
Let $A$ be a graded connected Gorenstein domain and also a factorial domain.
Then ${\rm Hom}_\F(A,\F)\cong A\theta_A$ with $\chi_A^{-1}=\lambda_A$ as defined
in Remark \ref{lambda_A}. Moreover $\Hom_\F(R,\F)\cong A_\lambda$, where
\begin{enumerate}
\item $\lambda:=\lambda_A$ if $\chr(k)$ does not divide $|G|$;
\item $\lambda:=\lambda_S\in {\rm Hom}(G/W,U(k))$, if $S=A^W$ is Cohen-Macaulay
(and therefore Gorenstein).
\end{enumerate}
In each of those cases $R=A^G$ is quasi-Gorenstein if and only if $\lambda=1$.
\end{prp}

\pf
It follows from \ref{main_intro} that there exists some function $\theta:=\theta_A$ with
$\Hom_\F(A,\F)\cong A\cdot \theta.$ From \cite{BH} Prop. 3.3.3 (a) we get
$$\overline A\cdot \overline{\theta}=
\Hom_\F(A,\F)\otimes \F/\F^+\cong \Hom_U(U,k)= U\cdot \xi_\ell^*,$$
hence $\overline{\theta}=c\cdot \xi_\ell^*$ with some nonzero scalar $c\in k$.
Setting $\lambda:=\lambda_A$, it follows that $\overline{g\theta}=g(\overline{\theta})=\lambda(g)^{-1}\overline\theta$,
so $g(\theta) - \lambda(g)^{-1}\cdot\theta\in \F^+\Hom_\F(A,\F).$
On the other hand $G$ maps $\theta$ onto another module generator and therefore
$g(\theta)=s_g\cdot\theta$ with a unit $s_g\in k=A_0$.
It follows that
$g(\theta) - \lambda(g)^{-1}\cdot\theta\in k\cdot\theta\cap \F^+A\theta=0$
and we conclude $g(\theta)=\lambda(g)^{-1}\cdot\theta$.
This shows $\chi_A=\lambda_A^{-1}$.\\
Since $S$ is a factorial domain, it is Gorenstein if Cohen-Macaulay, so the same argument
as above gives $\chi_S=\lambda_S^{-1}$. The statement about ${\rm Hom}_\F(R,\F)$
follows from \ref{main_intro}.\\
For the rest of the proof we can assume that $\chr(k)$ does not divide $|G|$.
We consider the restriction map
${\rm res}:\ \Hom_\F(A,\F)\to \Hom_\F(R,\F),\ \Psi\mapsto \Psi_{|R}$.
Since $t:\ A\to R,\ s\mapsto |G|^{-1}\sum_{g\in G} g(s)$ is an epimorphism
of $\F$-modules and $_\F R$ is free, we have
$_\F A =\ _\F R\oplus\ _\F X$ for some complement $_\F X\subseteq A$, hence
${\rm res}$ is surjective.
Let $\lambda\in \Hom_\F(R,\F)$, then there is $s\in A$ with
$\lambda=s\cdot \theta(|_R)=\theta(s\cdot())$. For any $r\in R$ we get
$$\lambda(r)=\frac{1}{|G|}\sum_{g\in G}\lambda(gr)=
\frac{1}{|G|}\sum_{g\in G}\theta(sgr)=
\frac{1}{|G|}\sum_{g\in G}\theta(g(g^{-1}(s)r))=$$
$$\frac{1}{|G|}\sum_{g\in G}\lambda(g)\theta(g^{-1}(s)r))=
\theta(t_\lambda(s)r),$$
where
$t_\lambda:=\frac{1}{|G|}\sum_{g\in G}\lambda(g)g^{-1}:\
A\to A_\lambda$
is the projection operator in $\Hom_\F(A,A_\lambda)$.
Thus we have $\Hom_\F(R,\F)\subseteq {\rm res}(A_\lambda\cdot \theta)$.
Again it follows from
$_\F A =\ _\F R\oplus\ _\F X$, that
$\Hom_\F(X,\F)^G=0$, hence
$$\Hom_\F(A,\F)^G\cong {\rm res}_R(\Hom_\F(A,\F)^G)=\Hom_\F(A^G,\F).$$
Clearly
$A_\lambda\cdot \theta\in \Hom_\F(A,\F)^G$, so
$$\Hom_\F(R,\F)= {\rm res}(A_\lambda\cdot \theta)\cong\ _\F A_\lambda.$$
If $\lambda=1$, then
$\Hom_\F(R,\F)= R\cdot {\rm res}(\theta)$ is
a cyclic $R$-module, so
$\omega_R\cong \Hom_\F(R,\F)\cong\ _RR$ and
$R$ is Gorenstein.
\epf

\begin{rem}
In the special case where $A={\rm Sym}(V^*)$ with linear $G$-action the
result above for the non-modular case also appears in
\cite{Nakajima:rel_inv} Cor. 3.2. The proof
indicated there depends on the results of
\cite{Watanabe:1}, \cite{Watanabe:2}. In contrast to this our proof above
is elementary and independent of Watanabe's results as well as of
our Theorem \ref{main_intro}.
\end{rem}

One can apply the results above for example in the situation where
$A:={\rm Sym}(V^*)$ for finite dimensional $kG$-module $V$, and
$S=A^W$ or $S=A^{\tilde W}$, with $\bar G:=G/W$ or $G/\tilde W$ acting
on $S$. However, even if $S=k[x_1,\cdots,x_n]$ is a polynomial ring (with ${\rm deg}(x_i)=:d_i\ge 1$),
then action of $\bar G$ will in general be non-linear and
the $k$-space $\langle x_1,\cdots,x_n\rangle_k$ will be not $\bar G$-stable.
Nevertheless we can use Remark \ref{Sym_V_example} to determine
$\lambda_S=\chi_S^{-1}$:

Let $M$ be a finite dimensional $kG$-module with $kG$-submodule $N\subseteq M$.
As a vectorspace we have $M=N\oplus U$, with $U\cong M/N$ as a $kG$-module.
Even though $M$ and $N\oplus U$ are in general not isomorphic as $kG$-modules,
one has $\chi_M=\chi_N+\chi_{M/N}$. It follows that
${\rm Sym}(M)\cong {\rm Sym}(N)\otimes_k {\rm Sym}(U)$ as a $k$-algebra, but in general
not as $kG$-module. Nevertheless we have
$H_{{\rm Sym}(M),g}^{Br}(t)=H_{{\rm Sym}(N),g}^{Br}(t)\cdot H_{{\rm Sym}(M/N),g}^{Br}(t).$
Even more generally, the following Lemma includes the case of a graded, but non-linear
$G$-action on the algebra generators:
\begin{lemma}\label{Brauer_series_of_tensor_prod}
Let $G$ act on $A$ by graded algebra automorphisms and $B\le A$ a $G$-stable graded subalgebra.
Assume that $A\cong B\otimes_k A/B_+A$ as a $k$-algebra (not necessarily as $kG$-module).
Then
$$H_{A,g}^{Br}(t)=H_{B,g}^{Br}(t)\cdot H_{A/B_+A,g}^{Br}(t).$$
\end{lemma}
\pf Let $A/B_+A=:C$ and identify the $k$-algebras
$B\otimes_k C\cong A,\ {\rm via}\  b\otimes c=bc.$
Let $x_1,\cdots,x_\mu$ be a $k$-basis of $B_m$ and
$y_1,\cdots,y_\nu$ a $k$-basis of $C_n$. Then
$g(y_j)=\sum_t g_{C;tj}y_t+\mathfrak{BC}$ with
$\mathfrak{BC}\in \sum_{r=1}^n B_rC_{n-r}$ and the matrix
$(g_{C;tj})$ describing the representation of $g$ on the $kG$-module
$C_j\cong (A/B_+A)_j$.
Hence $g(x_iy_j)=g(x_i)g(y_j)=$ $$\sum_s(g_{B;si}x_s)(\sum_t(g_{C;tj}y_t+\mathfrak{BC})=
g_{B;ii}\cdot g_{C;jj}\cdot x_iy_j+\sum_{(s,t)\ne(i,j)}g_{B;si}g_{C;tj}x_sy_t+ {\mathcal X},$$
with ${\mathcal X}\in \sum_{r=1}^n B_{m+r}C_{n-r}\subseteq B_+A.$
It follows that $\chi_{A_{m+n},g}=\chi_{B_m,g}\cdot \chi_{C_n,g}$ and therefore
$H_{A,g}^{Br}(t)=H_{B,g}^{Br}(t)\cdot H_{C,g}^{Br}(t).$
\epf

\begin{prp}\label{dets}
Let $A=k[x_{11},\cdots,x_{1j_1},x_{21},\cdots,x_{2j_2},\cdots,x_{\ell 1},\cdots,x_{\ell j_\ell}]$
be a polynomial ring with generators of degrees
$1\le d_1<d_2<\cdots<d_\ell.$ For $i:=1,\cdots,\ell$ let $U_i$ denote the $kG$-module 
$A_{d_i}/A^+A^+\cap A_{d_i}\in kG-mod$
and ${\rm det}_i:\ G\to k,\ g\mapsto \det(g_{|U_i})$.
Then for every $g\in G$:
$$H_{A,g}^{Br}(t)=\prod_{i=1}^k H_{{\rm Sym}(U_i),g}^{Br}(t)=
\prod_{i=1}^k \frac{1}{\widehat{\det(1-t^{d_i}g)}}\ {\rm and}$$
$\widehat{\lambda_A}(g)=\prod_{i=1}^k \widehat{\det}_i(g)^{-1}.$
\end{prp}
\pf The subalgebra
$B:=k[x_{11},\cdots,x_{1j_1}]={\rm Sym}(U_1)\subseteq A$ is $G$-stable
and we have $A=B\otimes_k A/B^+A$ with polynomial ring
$A/B^+A \cong k[\overline x_{21},\cdots,\overline x_{2j_2},\cdots,\overline x_{\ell 1},\cdots,\overline x_{\ell j_\ell}].$
Now the first equality follows from
Lemma \ref{Brauer_series_of_tensor_prod} and an obvious induction.
The rest follows in a way similar to Remark \ref{Sym_V_example}.
\epf

If $\chr(k)=p>0$, then by definition $p$ does not divide $[G:\tilde W]$, hence
if $A^{\tilde W}$ is Cohen-Macaulay, so is $A^G$. With regard to the Gorenstein
property we obtain the following:

\begin{cor}\label{w_tilde_cor}
Let $A:={\rm Sym}(V^*)$ with finite dimensional $kG$-module $V$ and
assume that $A^{\tilde W}\cong k[x_{11},\cdots,x_{1j_1},x_{21},\cdots,x_{2j_2},\cdots,x_{\ell 1},\cdots,x_{\ell j_\ell}]$
is a polynomial ring with generators of degrees
$1\le d_1<d_2<\cdots<d_\ell.$ Then $A^G$ is Cohen-Macaulay and $A^G$ is
Gorenstein if and only if $\prod_{i=1}^k \widehat{\det}_i(g)^{-1}=1$ for all $g\in G$.
\end{cor}

\bibliography{bibl}
\bibliographystyle{plain}
\end{document}